\documentclass[submission]{dmtcs}

\usepackage[english]{babel}
\usepackage{amsmath}
\usepackage{amsfonts}

\def\hs{\hbox to 3mm{}}
\def\hhs{\hbox to 5cm{}}

\author{K. A. Penson\addressmark{1}\thanks{Corresponding author: e-mail:penson@lptl.jussieu.fr}, P. Blasiak\addressmark{2},
G. H. E. Duchamp\addressmark{3}, A. Horzela\addressmark{2} and
A. I. Solomon\addressmark{1,4}.%\thanks{Corresponding author: K. A. Penson, e-mail:penson@lptl.jussieu.fr}
}
\title{On certain non-unique solutions of the Stieltjes moment problem}
\address{
\addressmark{1} Universit\'{e} Pierre et Marie Curie,
Laboratoire  de  Physique   Th\'{e}orique  de la Mati\`ere Condens\'ee, CNRS UMR 7600,
Tour 24, $2^{i\grave{e}me}$ \'{e}t., 4, Place Jussieu, F 75252 Paris Cedex 05, France.\\
\addressmark{2} H. Niewodnicza{\'n}ski Institute of Nuclear Physics,
Polish Academy of Sciences,
Division of Theoretical Physics,
ul. Eliasza-Radzikowskiego 152, PL 31-342 Krak{\'o}w, Poland.\\
\addressmark{3}LIPN, CNRS UMR 7030,
Institut Galil\'ee - Universit\'e Paris-Nord, 99, Av. Jean-Baptiste Cl\'ement, F 93430 Villetaneuse, France.\\
\addressmark{4}The Open University, Physics and Astronomy Department,
Milton Keynes MK7 6AA, United Kingdom.
}
\revision{26}
\received{\today}
\revised{\today}
\begin{document}
\maketitle
\date{}
\begin{abstract}
We construct explicit solutions of a number of Stieltjes moment problems based on moments of the form ${\rho}_{1}^{(r)}(n)=(2rn)!$ and ${\rho}_{2}^{(r)}(n)=[(rn)!]^{2}$, $r=1,2,\dots$, $n=0,1,2,\dots$, \textit{i.e.} we find functions $W^{(r)}_{1,2}(x)>0$ satisfying $\int_{0}^{\infty}x^{n}W^{(r)}_{1,2}(x)dx = {\rho}_{1,2}^{(r)}(n)$. It is shown using criteria for uniqueness and non-uniqueness (Carleman, Krein, Berg, Pakes, Stoyanov) that for $r>1$ both ${\rho}_{1,2}^{(r)}(n)$ give rise to non-unique solutions. Examples of such solutions are constructed using the technique of the inverse Mellin transform supplemented by a Mellin convolution. We outline a general method of generating non-unique solutions for moment problems generalizing ${\rho}_{1,2}^{(r)}(n)$, such as the product ${\rho}_{1}^{(r)}(n)\cdot{\rho}_{2}^{(r)}(n)$ and $[(rn)!]^{p}$, $p=3,4,\dots$.
\end{abstract}

\section{Introduction}\label{Introduction}
This paper concerns  solutions of the Stieltjes moment problem \cite{Akhiezer,Simon}, \textit{i.e.} positive functions $W(x)$ which satisfy the infinite set of equations
\begin{equation}\label{1}
\int\limits_{0}^{\infty}x^{n}W(x)dx = {\rho}(n), \ \ \ \ \ n=0,1,2,\dots .
\end{equation}
We previously met this  problem when considering the properties of the so-called {\em coherent states} (CS) of Quantum Mechanics \cite{cohstates}. These CS, as well as their many generalizations, should satisfy the property of {\em  Resolution of Unity} which is essentially equivalent to Eq. (\ref{1}) \cite{Klauder1,Klauder2}. Standard CS lead to $\rho(n)=n!$ and $W(x)=e^{-x}$. Generalized CS lead to $\rho(n)$'s other than $n!$ \cite{Barut} and the solutions of Eq.(\ref{1}), if any, must be studied in each individual case separately \cite{Penson,Penson1,Klauder3,Penson2}. It quickly occurred to us that an efficient approach to this problem was to use the inverse Mellin transform method which allows one to establish many solutions of Eq. (\ref{1}), either  by analytic methods \cite{Marichev} or by extensive use of available tables \cite{Prudnikov,Oberhettinger}. As a byproduct of this method we have established that, for a large number of combinatorial sequences such as  Bell and Catalan numbers,  \textit{etc.}, the corresponding sequences $\rho(n)$ are solutions of the moment problem Eq. (\ref{1}) \cite{Penson0, Penson01}. Likewise, sequences arising in theory of ordering of differential operators \cite{Penson4} solve appropriate moment problems too. This admixture of quantum-mechanical, analytical and combinatorial features deserves a deeper study which we intend to pursue.

Our paper is partly expository in character, and has the following structure:  first we establish a link between the Mellin transform and the moment problem;  next, in Section 3, we provide principal solutions to two moment problems, termed toy models. In Section 4 we discuss criteria for the uniqueness of solutions of the moment problem. Section 5 is devoted to the explicit construction of non-unique solutions of the toy-models. In Section 5 some generalizations of model sequences, together  with their solutions, are reviewed. Section 6 is devoted to a discussion and conclusions.
 
 {\em We dedicate this paper to Philippe Flajolet on the occasion of his 60th birthday. His pioneering applications of Mellin transform asymptotics  to the analysis of combinatorial structures \cite{Flajolet1,Flajolet2, Flajolet3} have been a  source of inspiration for us.}

\section{The Mellin transform \textit{versus} moment problem}
The Mellin transform of a function $f(x)$ of the real variable $x$ is defined for complex $s$ by the following relation
\begin{equation}
\label{m1}
{\cal M}[f(x);s] = \int\limits_{0}^{\infty}x^{s-1}f(x)dx\equiv f^{*}(s)
\end{equation}
( in this definition $f^{*}(s)$ is not  the complex conjugate of $f(s)$ ! ), and its inverse ${\cal M}^{-1}$ is defined by
\begin{equation}
\label{m2}
{\cal M}^{-1}[f^{*}(s);x] = \frac{1}{2\pi i}\int\limits_{c-i\infty}^{c+i\infty}f^{*}(s)x^{-s}ds.
\end{equation}
See Ref. \cite{Sneddon} for a discussion of the dependence of $f(x)$ on the real constant $c$. Among the many relations satisfied by the Mellin transform we shall mainly use the following
\begin{equation}
\label{m3}
{\cal M}[x^{b}f(ax^{h});s] = \frac{1}{h}a^{-\frac{s+b}{h}}f^{*}\left(\frac{s+b}{h}\right),\ \ \ \ a,h>0.
\end{equation}
If ${\cal M}[f(x);s]=f^{*}(s)$ and ${\cal M}[g(x);s]=g^{*}(s)$, then
\begin{equation}
\label{m4}
{\cal M}^{-1}[f^{*}(s)g^{*}(s);x] = \int\limits_{0}^{\infty}f\left(\frac{x}{t}\right)g(t)\frac{dt}{t} = \int\limits_{0}^{\infty}g\left(\frac{x}{t}\right)f(t)\frac{dt}{t}
\end{equation}
which is called the Mellin convolution property. Note that if in Eq. (\ref{m4}) both $f(x)$ and $g(x)$ are positive for $x>0$ then $\int\limits_{0}^{\infty}f(\frac{x}{t})g(t)\frac{dt}{t}$ is also positive for $x>0$. This means that the Mellin convolution preserves positivity, an essential property when considering the moment problem.
 
Eq.(\ref{1}), if rewritten for $n=s-1$ as
\begin{equation}
\label{m5}
W(x) = {\cal M}^{-1}[\rho(s-1);x],
\end{equation}
is, if $W(x)>0$, a solution of a Stieltjes moment problem. Thus according to Eq. (\ref{m5}) one can solve the Stieltjes moment problem by performing the inverse Mellin transform on the moment sequence and checking if the resulting function is positive. All the  solutions in the sequel have been obtained using Eqs. (\ref{m5}), (\ref{m3}) and (\ref{m4}). Note that $W(x)$ obtained \textit{via} Eq. (\ref{m5}) may not be the only solution of Eq.(\ref{1}). We call $W(x)>0$ obtained \textit{via} Eq. (\ref{m5}) from Eq. (\ref{1}) the principal solution.
\section{Principal solutions of the moment problems}\label{Principal}
Let us consider two sequences of integers given by
\begin{equation}
\label{2}
{\rho}_{1}^{(r)}(n)=(2rn)!,\ \ \ \ n=0,1,2,\dots,\ \ \ \ r=1,2,\dots.
\end{equation}
\begin{equation}
\label{3}
{\rho}_{2}^{(r)}(n)=[(rn)!]^{2},\ \ \ \ n=0,1,2,\dots,\ \ \ \ r=1,2,\dots.
\end{equation}
In the following we shall obtain the solutions of the Stieltjes moment problem for the moment sequences given by Eqs. (\ref{2}) and (\ref{3}), \textit{i.e.} the functions $W^{(r)}_{1,2}(x)>0$ satisfying
\begin{equation}
\label{4}
\int\limits^{\infty}_{0}x^{n}W^{(r)}_{1,2}(x)dx = \rho^{(r)}_{1,2}(n).
\end{equation}
From now on we shall refer to the model problems ${\rho}_{1,2}^{(r)}(n)$ as {\em toy models} \textbf{TM1} and \textbf{TM2}, respectively.

a) \textbf{TM1}: we begin by obtaining $W^{(1)}_{1}(x)$. Note that $(2n)!=\Gamma(2n+1)=\Gamma(2(s-\frac{1}{2}))$. We now apply  Eq.(\ref{m3}) with $a=1$, $b=-\frac{1}{2}$ and $h=\frac{1}{2}$. We subsequently use ${\cal M}^{-1}[\Gamma(s);x] = e^{-x}$ which gives
\begin{equation}
\label{5}
\int\limits^{\infty}_{0}x^{n}W^{(1)}_{1}(x)dx = \int\limits^{\infty}_{0}x^{n}\left[\frac{e^{-\sqrt{x}}}{2\sqrt{x}}\right]dx = (2n)!,\ \ \ \ n=0,1,2,\dots.
\end{equation}
In the same spirit we observe that $(2rn)! = \Gamma\left(2r\left(n+\frac{1}{2r}\right)\right)=\Gamma\left(2r\left(s-\frac{2r-1}{2r}\right)\right)$. Upon using Eq. (\ref{m3}) but now with $a=1$, $b=-\frac{2r-1}{2r}$ and $h=\frac{1}{2r}$ one obtains
\begin{equation}
\label{6}
\int\limits^{\infty}_{0}x^{n}W^{(r)}_{1}(x)dx = \int\limits^{\infty}_{0}x^{n}\left[\frac{1}{2rx^{\frac{2r-1}{2r}}}e^{-x^{\frac{1}{2r}}}\right]dx = (2rn)!,\ \ \ \ n=0,1,2,\dots.
\end{equation}
We call $W^{(r)}_{1}(x)>0$ the {\em principal solution} of the moment problem Eq. (\ref{6}).

b) \textbf{TM2}: we begin by deriving $W^{(1)}_{2}(x)={\cal M}^{-1}[\Gamma^{2}(s);x]$ and employ the Mellin convolution Eq. (\ref{m4}). By using the Sommerfeld representation of the modified Bessel function of second kind $K_{0}(x)$ \cite{Lebedev} one obtains
\begin{equation}
\label{7}
\int\limits^{\infty}_{0}x^{n}W^{(1)}_{2}(x)dx = \int\limits^{\infty}_{0}x^{n}\left[2K_{0}(2x^{\frac{1}{2}})\right]dx = (n!)^2,\ \ \ \ n=0,1,2,\dots.
\end{equation}
Subsequently note that $[(rn)!]^2 = [\Gamma\left(r\left(n+\frac{1}{r}\right)\right)]^2 = [\Gamma\left(r\left(s-\frac{r-1}{r}\right)\right)]^2$ and  again apply Eq.(\ref{m3}) with $a=1$, $b=-\frac{r-1}{r}$ and $h=\frac{1}{r}$. The result is
\begin{equation}
\label{8}
\int\limits^{\infty}_{0}x^{n}W^{(r)}_{2}(x)dx = \int\limits^{\infty}_{0}x^{n}\left[\frac{2}{rx^{\frac{r-1}{r}}}K_{0}(2x^{\frac{1}{2r}})\right]dx = [(rn)!]^2,\ \ \ \ n=0,1,2,\dots.
\end{equation}
Since $K_{0}(t)>0$ for $t>0$, $W_{2}^{(r)}(x)>0$ is the principal solution of the moment problem Eq. (\ref{8}).
\section{Criteria of uniqueness and non-uniqueness of the Stieltjes moment problem}
It was realized from the very beginning of the history of the moment problem that its solutions may not be unique; \textit{i.e.} for a given moment sequence there may exist more than one solution. Stieltjes himself gave an example of a non-unique solution of a problem leading to what turned out to be the lognormal distribution \cite{Stieltjes}. This example, being  about the only one available,  was quoted repeatedly in the literature. As recently as  twenty years ago new non-unique solutions of other types of problems have been constructed. Of special interest were Stieltjes moment problems arising in probability theory, related to investigation of probability distributions not determined by their moments. Consequently, the subject was further developed and systematized, largely due to the comprehensive work of Berg \cite{Berg}, Berg and Pedersen \cite{Pedersen}, Lin \cite{Lin}, Stoyanov \cite{Stoyanov2,Stoyanov1}, Pakes \cite{Pakes}, Gut \cite{Gut} and others. From the practical point of view one needs criteria to decide whether the pursuit of non-unique solutions is reasonable. Such criteria are either based on the $\rho(n)$'s alone or on the solution $W(x)$ alone, or  on both $\rho(n)$ and $W(x)$, see below. From now on we assume that all the $\rho(n)$'s are finite and that $W(x)$ is continuous.
We now give , in a somewhat condensed form, a list of such criteria. 
\begin{description}
\item\textbf{C1} Carleman uniqueness criterion (T. Carleman, 1922, \cite{Akhiezer})
This  is based on the properties of the $\rho(n)$'s alone and is:

{\it If $S=\sum_{n=1}^{\infty}[\rho(n)]^{-\frac{1}{2n}} = \infty$, then the solution is unique.}

This criterion does not imply that if $S<\infty$ then the solution is non-unique. In fact it is possible to construct models for which $S<\infty$ and solutions are still unique \cite{Stoyanov2,Simon1}.

\item\textbf{C2} Krein's non-uniqueness criterion (M. G. Krein, around 1950, \cite{Krein})

 This is  based entirely on the solution $W(x)$ and does not involve the moments.

{\it If $\int_{0}^{\infty}\frac{-\ln{[W(x^2)]}}{1 + x^2}dx < \infty$, \textit{i.e.}, the so-called Krein integral exists, then the solution is non-unique.}

\item\textbf{C3} Converse Carleman criterion for non-uniqueness (A. Pakes, 2001, \cite{Pakes}, A. Gut, 2002, \cite{Gut})

This is based on $\rho(n)$'s and $W(x)$.

{\it If there exists $x^{\prime}\ge 0$ such that for $x>x^{\prime}$, $0<W(x)<\infty$ and $\psi(y)=-\ln{[W(e^{y})]}$ is convex in $(y^{\prime},\infty)$, where $y^{\prime}=\ln{(x^{\prime})}$, and if, in addition, $S=\sum_{n=1}^{\infty}[\rho(n)]^{-\frac{1}{2n}} < \infty$, then the solution $W(x)$ is non-unique.}
\end{description}
See \cite{Berg,Pedersen,Stoyanov2,Pakes,Gut} for various refinements of these criteria.

\section{Construction of non-unique solutions for TM1 and TM2}
We first apply the criterion C1 to sequences $\rho_{1,2}^{(r)}(n)$ (the logarithmic test of divergence of the series $S$ is conclusive) and conclude that for $r=1$ both solutions $W^{(1)}_{1}(x)=\frac{e^{-\sqrt{x}}}{2\sqrt{x}}$ and $W^{(1)}_{2}(x) = 2K_{0}(2x^{\frac{1}{2}})$ are unique. For $r>1$ we find that $\sum_{n=1}^{\infty}[\rho^{(r)}_{1,2}(n)]^{-\frac{1}{2n}}$ is convergent. Application of criterion C2 gives convergence of the Krein integral, showing that the solutions of \textbf{TM1} and \textbf{TM2} are non-unique. These findings are confirmed by use of criterion C3: convexity of $\psi_{1,2}^{(r)}(x)=-\ln{[W^{(r)}_{1,2}(e^x)]}$ is proved as it is equivalent to $e^{x/r}>0$ for $x>0$, (\textbf{TM1}) and $K_{1}(2e^{x/r})-K_{0}(2e^{x/r})>0$ for $x>0$, (\textbf{TM2}). The quest for non-unique solutions for $r>1$ is then well founded.

We know of no general method to construct such solutions. However we have proposed a procedure based on the application of inverse Mellin transform which generates the required solutions \cite{Penson1,Sixdeniers} which we now briefly expose.

The first step is to construct, within the framework of a given set of $\rho(n)$'s, a family of functions $\omega_{k}(x)$, parametrized by a constant $k$ (to be defined below), such that \textit{all} their moments vanish, \textit{i.e.}
\begin{equation}
\label{9}
\int\limits^{\infty}_{0}x^{n}\omega_{k}(x)dx = \int\limits^{\infty}_{0}x^{s-1}\omega_{k}(x)dx = 0, \ \ \ \ n=0,1,2,\dots,\ \ \ \ s=1,2,\dots.
\end{equation}
The Mellin transform of $\omega_{k}(x)$, $\omega^{\star}_{k}(s)$, vanishes for $s=1,2,\dots$. Such functions are orthogonal to all polynomials and play an important role in the study of integral transforms \cite{Samko}. For our purposes we choose a particular method of producing the functions $\omega^{\star}_{k}(x)$:
\begin{equation}
\label{10}
\int\limits^{\infty}_{0}x^{n}\omega_{k}(x)dx = \rho_{1,2}^{(r)}(n)\cdot h_{k}(n),
\end{equation}
or equivalently
\begin{equation}
\label{11}
\int\limits^{\infty}_{0}x^{s-1}\omega_{k}(x)dx = \rho_{1,2}^{(r)}(s-1)\cdot h_{k}(s-1),
\end{equation}
where $h_{k}(s)$ is any holomorphic function vanishing for $s=1,2,\dots$. Among an infinity of possible choices the simplest one is $h_{k}(s)\!=\!\sin{(\pi k(s+1))}$ and it defines a discrete parameter $k\!=\!\pm 1,\pm 2, \pm 3,\dots$. The function $\omega_{k}(x)$ acquires new parameters now and is formally obtained by calculating the inverse Mellin transform
\begin{equation}
\label{12}
\omega_{1,2,k}^{(r)}(x)=\frac{1}{2\pi i}\int\limits^{i\infty}_{-i\infty}\rho_{1,2}^{(r)}(s-1)\sin{(\pi ks)}x^{-s}ds,\ \ \ \ k=\pm 1, \pm 2, \pm 3, \dots.
\end{equation}
It turns our that for both $\rho_{1,2}^{(r)}(s)$ the integration in Eq.(\ref{12}) can be performed:

\noindent a) \textbf{TM1}: for $\rho_{1}^{(r)}(n)=(2rn)!$, $\omega^{(r)}_{1,k}(x)$ is a special case of earlier evaluation \cite{Penson1} and reads:
\begin{equation}
\label{13}
\begin{array}{rcl}
\omega^{(r)}_{1,k}(x) &=& \frac{1}{2rx^{\frac{2r-1}{2r}}}e^{-x^{\frac{1}{2r}}}\sin{\left[k\pi\left(\frac{2r-1}{2r}\right) + x^{1/2r}\tan{\left(\frac{k\pi}{2r}\right)}\right]}, \ \ \ \ (r>|k|)
\\\\
&=&W_{1}^{(r)}(x)\sin{\left[k\pi\left(\frac{2r-1}{2r}\right) + x^{1/2r}\tan{\left(\frac{k\pi}{2r}\right)}\right]}.
\end{array}
\end{equation}
In Eq.(\ref{13}) we notice a pleasant factorization of $W_{1}^{(r)}(x)$.

\noindent b) \textbf{TM2}: for $\rho_{2}^{(r)}(n)= [(rn)!]^2$ the corresponding function $\omega^{(r)}_{2,k}(x)$ has to be, in the first place, represented as
\begin{equation}
\label{14}
\omega^{(r)}_{2,k}(x)={\cal{M}}^{-1}[\underbrace{\Gamma(rs-s+1)}_{\textrm{I}}\underbrace{\Gamma(rs -s+1)\sin{(\pi ks)}}_{\textrm {II}};x]
\end{equation}
which can be conceived as another case of Mellin convolution. A little thought gives the two partners to be convoluted as
\begin{equation}
\label{15}
\textrm{I}\rightarrow W_{1}^{(r/2)}(x) = \frac{1}{rx^{\frac{r-1}{r}}}e^{-x^{\frac{1}{r}}},
\end{equation}
(see Eq.(\ref{13})) and
\begin{equation}
\label{16}
\begin{array}{c}
\textrm{II}\rightarrow \omega_{1,k}^{(r/2)}(x) = \frac{1}{rx^{\frac{r-1}{r}}}e^{-x^{\frac{1}{r}}}\sin{\left[k\pi\left(\frac{r-1}{r}\right) + x^{1/r}\tan{\left(\frac{k\pi}{r}\right)}\right]},
\end{array}
\end{equation}
(compare Eq.(\ref{13})). Thus, the integral form of  $\omega_{2,k}^{(r)}(x)$ is
\begin{equation}
\label{17}
\omega_{2,k}^{(r)}(x) = \int\limits_{0}^{\infty}W_{1}^{(r/2)}\left(\frac{x}{t}\right)\omega_{1,k}^{(r/2)}(t)\frac{dt}{t}
\end{equation}
whose evaluation requires a number of changes of variables as well as the use of formula 2.5.37.2, p. 453 of vol.1 Ref.\cite{Prudnikov}, but is essentially elementary. The final result is
\begin{equation}
\label{18}
\begin{array}{rcl}
\omega_{2,k}^{(r)}(x) &=& \frac{2}{rx^{\frac{r-1}{r}}}\textrm{Re}\left[e^{i\pi\left(\frac{1}{2}-k\frac{r-1}{r}\right)}
K_{0}\!\left(2x^{1/2r}\left(1 + i\tan\left(\frac{\pi k}{r}\right)\right)^{1/2}\right)\right]
\\\\
&\equiv& \frac{2}{rx^{\frac{r-1}{r}}}V^{(r)}_{k}(x)
\end{array}
\end{equation}
where we note a ``near'' factorization of $W^{(r)}_{2}(x)$.

Armed with explicit forms for $\omega_{1,k}^{(r)}(x)$ and $\omega_{2,k}^{(r)}(x)$ we are in position now to write down the families of non-unique solutions. Their structure has the form: principal solution + const$\cdot\omega_{k}(x)$. More precisely:

\noindent\textbf{TM1}: 
\begin{equation}
\label{18a}
{\tilde W}^{(r)}_{1}(\epsilon,k,x) = W^{(r)}_{1}(x)\left[1 + \epsilon\sin{\left(k\pi\left(\frac{2r-1}{2r}\right)+x^{1/2r}\tan{\left(\frac{k\pi}{2r}\right)}\right)}\right]
\end{equation} 
for real $\epsilon$, $|\epsilon|<1$.

\noindent\textbf{TM2}: 
\begin{equation}
\label{18b}
{\tilde W}^{(r)}_{2}(\gamma,k,x) = W^{(r)}_{2}(x)\left[1 + \gamma\frac{V_{k}^{(r)}(x)}{K_{0}(2x^{1/2r})}\right]
\end{equation}
for $r>2|k|$.

\noindent As $\left[1 + \gamma\frac{V_{k}^{(r)}(x)}{K_{0}(2x^{1/2r})}\right]$ is an oscillating function of bounded variation, a constant $\gamma=\gamma(k,r)$ can be always found to assure the overall positivity of ${\tilde W}^{(r)}_{2}(\gamma,k,x)$. The above technique for obtaining non-unique solutions can be readily extended to moment sequences more general than $\rho_{1,2}^{(r)}(n)$. We shall simply mention two such extensions without entering into details. 

For $\rho_{3}^{(r)}(n)=[(rn)!]^{3}$ we begin with the sequence $(n!)^{3}$ for which the solution is
\begin{equation}
\label{19}
\int\limits_{0}^{\infty}x^{n}\textrm{MeijerG}\left([\ [\ \ ],[\ \ ]\ ],[\ [0,0,0],[\ \ ]\ ],x\right)dx = (n!)^3,
\end{equation}
where we use a convenient and self-explanatory notation for Meijer's G-function borrowed from that of computer algebra systems. The extension, \textit{via} Eq. (\ref{m3}), leads to the principal solution
\begin{equation}
\label{20}
\int\limits_{0}^{\infty}x^{n}\left[\frac{1}{rx^{\frac{r-1}{r}}}\textrm{MeijerG}\left([\ [\ \ ],[\ \ ]\ ],[[0,0,0],[\ \ ]],x^{1/r}\right)\right]dx = [(rn)!]^3.
\end{equation}
The integrand in Eq. (\ref{20}) is a positive function for $x>0$ which cannot be represented by any other known special function. It possesses an infinite series representation in terms of polygamma functions, which we will not quote here. The Carleman sum $S$ is convergent but the criterion C2 is not conclusive. Only the criterion C3 permits to ascertain the non-uniqueness. The corresponding function  $\omega_{3,k}^{(r)}(x)$ is defined as
\begin{equation}
\label{21}
\omega_{3,k}^{(r)}(x) = {\cal M}^{-1}[\underbrace{\Gamma^{2}(rs - s +1)}_{\textrm{I}}\underbrace{\Gamma(rs - s +1)\sin{(\pi ks)}}_{\textrm{II}};x]
\end{equation}
which can be calculated as Mellin convolution of $\textrm{I}\rightarrow W^{(r)}_{2}(x)$ and $\textrm{II}\rightarrow \omega_{1,k}^{(r/2)}(x)$, see Eqs. (\ref{13}) and (\ref{18}), respectively.

As a final example consider the sequence $\rho_{4}^{(r)}(n)=\rho_{1}^{(r)}(n)\rho_{2}^{(r)}(n)$. The corresponding Mellin convolution of two principal solutions $W^{(r)}_{1,2}(x)$, see Eqs. (\ref{18a}) and (\ref{18b}), yields directly the principal solution for $\rho_{4}^{(r)}(n)$:
\begin{equation}
\label{22}
\begin{array}{rl}
\int\limits_{0}^{\infty}x^{n}\left[\frac{4^{r-1}}{r\sqrt{\pi}x^{\frac{2r-1}{r}}}\textrm{MeijerG}\left([\ [\ \ ],[\ \ ]\ ],[[r-\frac{1}{2},r,r,r],[\ \ ]],\frac{1}{4}x^{1/r}\right)\right]dx &= (2rn)![(rn)!]^2,
\\  
&n=0,1,2\dots,\ \ r=1,2,\dots,
\end{array}
\end{equation}
which is non-unique by C3 only, as C2 remains inconlusive.
\section{Discussion and Conclusion}
We have demonstrated a methodology for obtaining unique and non-unique solutions of the Stieltjes moment problem using the Mellin convolution method. Although the initial moment sequences were simple and classical, one is rapidly forced to leave the realm of standard special functions, as the resulting solutions are special cases of Meijer G-functions. In most cases they resist the check for non-uniqueness \textit{via} both the Carleman criterion C1 and the Krein criterion C2, and criterion C3 appears to be the only tool to decide this question. We have generated parametrized families of non-unique solutions exemplified here by Eqs. (\ref{13}) and (\ref{18}). Such functions are now called \textit{Stieltjes classes} \cite{Stoyanov3, Pakes2} and their properties are investigated.

\section*{Acknowledgments}
We thank A. Gut, B. Chemin and H. L. Pedersen for discussions. We also wish to acknowledge support from Agence Nationale de la Recherche (Paris, France) under Program No. ANR-08-BLAN-0243-2 and from PAN/CNRS Project PICS No.4339 (2008-2010). Two of us (P.B. and A.H.) wish to acknowledge support from Polish Ministry of Science and Higher Education under Grants Nos. N202 061434 and N202 107 32/2832.


\begin{thebibliography}{0}
\bibitem{Krein} M. G. Krein formulated logarithmic integral criteria for the Hamburger moment problem in 40' and 50' of the XXth century. Extensions of Krein's criteria to the Stieltjes moment problem were found much later.
\end{thebibliography}
\end{document}